\documentclass{article}

\usepackage{amsmath,amssymb,amsthm}
\usepackage{latexsym}
\usepackage{enumerate}
\usepackage{geometry}
\usepackage{hyperref}
\usepackage{cite}

\theoremstyle{plain}

\theoremstyle{definition}

\theoremstyle{remark}

\newcommand{\K}[1]{\textit{K}_{#1}}
\newcommand{\B}[1]{\textit{B}_{#1}}
\newcommand{\W}[1]{\textit{W}_{#1}}
\newcommand{\BP}{\textbf{P}}
\newcommand{\Ag}{\textbf{I}}

\newcommand{\M}{\mathcal{M}}
\newcommand{\forces}{\vDash}
\newcommand{\nforces}{\nvDash}
\newcommand{\all}{\forall}
\newcommand{\some}{\exists}
\newcommand{\Dia}{\Diamond}
\newcommand{\Boxm}{\Box}

\title{\textbf{Taming Fitch's and Humberstone's Paradoxes}}
\author{Jie Fan\\
\normalsize Institute of Philosophy, Chinese Academy of Sciences;\\
\normalsize School of Humanities, University of Chinese Academy of Sciences\\
\normalsize \texttt{jiefan@ucas.ac.cn}}
\date{}

\begin{document}
\maketitle

\begin{abstract}
Fitch's Paradox says that if every true proposition is knowable, then every true proposition is known. Humberstone's Paradox states that if every true proposition is known by at least one of a set of finitely many agents, then some agent in that set knows every true proposition. Both paradoxes give rise to surprising and unacceptable consequences. In this paper, we diagnose a crucial tacit presupposition underlying these paradoxes, namely that the domain of propositions is closed under Boolean operations such as intersection and union, and offer a uniform resolution. We argue that the paradoxical inferences rely implicitly on this closure assumption. By adopting weak propositionally quantified modal logics that drop these closure conditions, we provide a semantically motivated account that blocks both paradoxes. The restriction on universal instantiation is not an ad hoc technical device, but reflects a conception of propositions as graspable objects of epistemic attitudes, on which the domain of propositions need not be closed under Boolean operations, since some Boolean combinations of graspable propositions may not themselves be graspable.
\end{abstract}

\noindent\textbf{Keywords:} Fitch's Paradox, Humberstone's Paradox, propositional quantified modal logic

\section{Introduction}

In a seminal work \cite{fitch1963}, Fitch shows that if every true proposition is knowable, then every true proposition is known; or equivalently, if there is an unknown truth, then there is an unknowable truth. This is known in the literature as ``Fitch's Paradox'' (also called ``Fitch's Paradox of Knowability'', ``the Knowability Paradox'', or ``Church-Fitch Paradox''), see e.g. \cite{brogaard2019}. In a parallel strand of the literature, Humberstone \cite{humberstone1985} argues that if every true proposition is known by at least one of a set of finitely many agents, then some agent in that set knows every true proposition. This is called ``The Humberstone Paradox'' \cite{bigelow2005} or ``Humberstone's Paradox'' \cite{maglo2010,updike2024}.

Both paradoxes give rise to surprising and unacceptable consequences. The knowability principle, namely, that every true proposition is knowable, has been claimed for a multitude of historical non-realisms, for instance Michael Dummett's semantic anti-realism, Hilary Putnam's internal realism, the logical positivisms of the Berlin and Vienna Circles, Peirce's pragmatism, Kant's transcendental idealism, and Berkeley's metaphysical idealism. Fitch's Paradox takes us from the ``relatively moderate and plausible'' knowability principle to the ``exceedingly implausible and obviously false omniscience principle'' (namely, that every true proposition is actually known) \cite[pp. 1--2]{salerno2009book}, thereby collapsing the distinction between truth and knowledge. This is untenable, given that someone committed to the claim that every true proposition is knowable should not thereby automatically be committed to the claim that every true proposition is known. So we need to find a way to block Fitch's argument aiming to show that the first of these claims entails the second. Likewise, it is plausible that there is a distribution of shared space-time resources whereby all agents collectively could come to know (via both individual and joint efforts) all true propositions without any single agent becoming actually omniscient \cite[pp. 1183--1184]{updike2024}, a scenario that contradicts Humberstone's Paradox.

Since Fitch's 1963 publication, a substantial literature has emerged addressing Fitch's paradox; see \cite{brogaard2019,salerno2009book} and the references therein. Many treatments have been proposed, including logical revisions (such as epistemic, intuitionistic, and paraconsistent revisions), semantic restrictions, and syntactic restrictions, with \cite{brogaard2019} providing a comprehensive survey.

Compared to the well-known Fitch's Paradox, Humberstone's Paradox has received remarkably little attention. Adapting Humberstone's argument, Bigelow \cite{bigelow2005} argues that there is a God.\footnote{But see Guigon \cite{guigon2009} and Humberstone \cite[Chapter 6]{humberstone2015book}, who challenge the assumption that the bring-about operator distributes over conjunction. Even so, Updike \cite{updike2024} recovers Bigelow's results in a framework where this operator does not distribute over conjunction.} Updike \cite{updike2024} argues that, even without the assumption that knowledge distributes over conjunction (namely, $\K{a}(p\land q)\to(\K{a}p\land\K{a}q)$), Humberstone's Paradox can still be recovered. Both works presuppose the validity of Humberstone's Paradox. This paradox is related to areas such as action theory and theology, see \cite{linsky1986,maglo2010,restall2009,updike2022,updike2024}. However, given its counterintuitive and unacceptable results (e.g. \cite{bigelow2005}), we have reason to seek to avoid the Paradox. A natural question thus arises: is there any approach to block Humberstone's Paradox?

While numerous solutions exist for Fitch's Paradox, many face the charge of being \emph{ad hoc}, lacking independent philosophical motivation beyond blocking the paradoxical inference, and Humberstone's Paradox remains largely unaddressed. In this paper, we offer a uniform resolution to both paradoxes, grounded in a minimal and plausible semantics of propositional quantification that explains why the paradoxical steps are illegitimate rather than merely declaring them invalid. More specifically, by reconstructing Humberstone's Paradox in propositionally quantified modal logic, we show that Humberstone's argument does not hold in general, since the initial transition is blocked. We also show that this approach extends to Fitch's Paradox and to the variations of Humberstone's Paradox. We conclude by offering a detailed response to potential objections concerning the philosophical motivation of our framework, clarifying why our restriction on universal instantiation constitutes a principled resolution rather than an \emph{ad hoc} technical maneuver.

The remainder of the paper is organized as follows. Section~\ref{sec.humberstoneparadox} introduces Humberstone's Paradox, its variations, and their original arguments. Section~\ref{sec.fitchparadox} introduces Fitch's Paradox and its original argument. In Section~\ref{sec.approach}, we present a uniform approach to resolving these paradoxes. Finally, Section~\ref{sec.conclusion} addresses a plausible objection and offers our responses.

\section{Humberstone's Paradox and Its Variations}\label{sec.humberstoneparadox}

Humberstone's argument, from collective knowledge to individual omniscience, can be expressed formally. For simplicity, we restrict attention to the two-agent case,\footnote{The argument extends straightforwardly to the case of finitely many agents.} which can be formulated in a bimodal logic as follows: if $p\to(\K{a}p\lor\K{b}p)$, then $(p\to\K{a}p)\lor(q\to\K{b}q)$, where $\K{a}$ and $\K{b}$ are the standard knowledge operators, meaning ``$a$ knows that'' and ``$b$ knows that'', respectively. Humberstone's original argument is as follows.

\begin{enumerate}[(1)]
\item $p\to(\K{a}p\lor\K{b}p)$ \hfill Assumption
\item $(p\land q)\to(\K{a}(p\land q)\lor\K{b}(p\land q))$ \hfill (1), substituting $p\land q$ for $p$
\item $\K{a}(p\land q)\to\K{a}p$ \hfill By Monotonicity applied to $(p\land q)\to p$
\item $\K{b}(p\land q)\to\K{b}q$ \hfill By Monotonicity applied to $(p\land q)\to q$
\item $(p\land q)\to(\K{a}p\lor\K{b}q)$ \hfill A truth-functional consequence of (2)--(4)
\item $(p\to\K{a}p)\lor(q\to\K{b}q)$ \hfill Truth-functionally equivalent to (5)
\end{enumerate}

Humberstone extends the argument to other contexts. For instance, by replacing $\K{}$ with the permissibility operator $P$ and reinterpreting $a$ ($b$) as a moral code rather than an agent, we obtain $P_a\varphi$ (resp. $P_b\varphi$), read as ``$\varphi$ is permissible according to moral code $a$ (resp. $b$).'' Since permissibility is plausibly monotonic, an argument analogous to (1)--(6) yields the following: if every true proposition is permissible according to one of the two codes or the other (in symbols, $p\to(P_ap\lor P_bp)$), then there is at least one code according to which every true proposition is permissible (in symbols, $(p\to P_ap)\lor(q\to P_bq)$).

As another example, by replacing $p$ in the antecedent of the above (1) with $\K{c}p$, we obtain $\K{c}p\to(\K{a}p\lor\K{b}p)$, which is read as ``$a$ and $b$ are collectively as knowledgeable as $c$''. A similar argument yields $(\K{c}p\to\K{a}p)\lor(\K{c}q\to\K{b}q)$, which says that either $a$ alone or $b$ alone is at least as knowledgeable as $c$.

If knowledge is replaced with belief in the above argument, then (1) says that $a$ and $b$ are collectively successful in having true beliefs, while (6) says that either $a$ or $b$ alone has achieved such success. Humberstone then considers the analogue of (1) for the goal of avoiding false beliefs, namely the assumption that whatever both believe is true. At first glance, this seems to allow each to have false beliefs not shared by the other. But surprisingly, this is not so: if whatever $a$ and $b$ both believe is true, then at least one of them has only true beliefs. The proof goes as follows, where for any proposition $\varphi$, $\B{i}\varphi$ is read as ``the agent $i$ believes that $\varphi$''.

\begin{enumerate}[(1$'$)]
\item $(\B{a}p\land\B{b}p)\to p$ \hfill Assumption
\item $(\B{a}(p\lor q)\land\B{b}(p\lor q))\to(p\lor q)$ \hfill (1$'$), substituting $p\lor q$ for $p$
\item $\B{a}p\to\B{a}(p\lor q)$ \hfill By Monotonicity applied to $p\to(p\lor q)$
\item $\B{b}q\to\B{b}(p\lor q)$ \hfill By Monotonicity applied to $q\to(p\lor q)$
\item $(\B{a}p\land\B{b}q)\to(p\lor q)$ \hfill A truth-functional consequence of (2$'$)--(4$'$)
\item $(\B{a}p\to p)\lor(\B{b}q\to q)$ \hfill Truth-functionally equivalent to (5$'$)
\end{enumerate}

Then, as Humberstone remarks, if we use $\W{i}\varphi$ to represent that $i$ wants it to be the case that $\varphi$, and assume monotonicity for the operator $\W{}$, then the above argument from (1$'$) to (6$'$), with $\B{}$ replaced by $\W{}$, gives us another surprising result: if whatever $a$ and $b$ both want to be the case is the case (in symbols, $(\W{a}p\land\W{b}p)\to p$), then either all $a$'s desires or all $b$'s desires are satisfied (in symbols, $(\W{a}p\to p)\lor(\W{b}q\to q)$). This entails that if a group of finitely many agents is collectively omnipotent (better termed ``omnicontented''), then at least one member of the group is omnipotent outright.

\section{Fitch's Paradox}\label{sec.fitchparadox}

Fitch's argument can be expressed in a modal epistemic logic: from $p\to\Dia\K{}p$ it follows that $p\to\K{}p$, where $\Dia$ means ``it is possible that'' and $\K{}$ means ``it is known that''. The argument proceeds as follows.\footnote{In fact, Fitch's original result is that if there is an unknown truth, then there is an unknowable truth; equivalently, knowability implies omniscience. Fitch credits the proof of \cite[Thm.~4]{fitch1963} to an anonymous referee, identified by Salerno \cite{salerno2009} as Alonzo Church. Proofs vary slightly in the literature. For instance, rather than deriving (7$''$) as we do here, van Benthem \cite{vanbenthem2009} uses the monotonicity of $\Dia$: $\Dia\K{}(p\land\neg\K{}p)\to\Dia(\K{}p\land\K{}\neg\K{}p)\to\Dia(\K{}p\land\neg\K{}p)\to\Dia\bot\to\bot$. This variant is inessential for our purposes; we adopt Dummett's formulation from \cite{dummett2009}.}

\begin{enumerate}[(1$''$)]
\item $p\to\Dia\K{}p$ \hfill Assumption
\item $(p\land\neg\K{}p)\to\Dia\K{}(p\land\neg\K{}p)$ \hfill (1$''$), substituting $p\land\neg\K{}p$ for $p$
\item $\K{}(p\land\neg\K{}p)\to(\K{}p\land\K{}\neg\K{}p)$ \hfill By distribution of $\K{}$ over $\land$
\item $\K{}\neg\K{}p\to\neg\K{}p$ \hfill Factivity of knowledge
\item $\K{}(p\land\neg\K{}p)\to(\K{}p\land\neg\K{}p)$ \hfill A truth-functional consequence of (3$''$) and (4$''$)
\item $\neg(\K{}p\land\neg\K{}p)$ \hfill Tautology
\item $\neg\K{}(p\land\neg\K{}p)$ \hfill A truth-functional consequence of (5$''$) and (6$''$)
\item $\Boxm\neg\K{}(p\land\neg\K{}p)$ \hfill (7$''$), Necessitation
\item $\neg\Dia\K{}(p\land\neg\K{}p)$ \hfill Equivalent to (8$''$)
\item $\neg(p\land\neg\K{}p)$ \hfill A truth-functional consequence of (9$''$) and (2$''$)
\item $p\to\K{}p$ \hfill Truth-functionally equivalent to (10$''$)
\end{enumerate}

Fitch \cite{fitch1963} and others (e.g. Salerno \cite{salerno2009}) observe that the reasoning in Fitch's argument generalizes to any operator that is factive and distributes over conjunction. Equivalently, for any operator $O$ satisfying $Op\to p$ and $O(p\land q)\to Op\land Oq$, if $p\to\Dia Op$, then $p\to Op$.\footnote{Note that neither distribution of $\K{}$ over $\land$ nor factivity is strictly required for the above derivation. A belief-version uses the KD4 theorem $\B{}\neg\B{}p\to\neg\B{}p$ in place of factivity; see e.g. \cite{linsky1986}, which employs a generalized version of this theorem. A knowledge-version also avoids distribution of $\K{}$ over $\land$; see e.g. \cite{maglo2010,williamson1993}. These variants of Fitch's Paradox can also be handled by our approach with suitable adaptations.}

\section{A Perspective on Propositional Quantification}\label{sec.approach}

The careful reader will notice that both Humberstone's and Fitch's arguments rely on implicit propositional quantifiers. When the quantifiers are made explicit, Humberstone's Paradox takes the form: if $\all p(p\to\K{a}p\lor\K{b}p)$, then $\all p(p\to\K{a}p)\lor\all q(q\to\K{b}q)$, and Fitch's Paradox takes the form: if $\all p(p\to\Dia\K{}p)$, then $\all p(p\to\K{}p)$. These implicit quantifiers motivate our uniform approach to blocking both paradoxes. In what follows, we first resolve Humberstone's Paradox and its variations, then Fitch's Paradox.

The notion of propositional quantification is introduced by Fine \cite{fine1970}, who provides strong axiomatizations. We do not need the full strength of these systems for our purposes; instead, we use a weak propositionally quantified modal logic containing only the axioms and inference rules needed below. Specifically, let $\mathbb{S}$ consist of the following axioms and inference rules:\footnote{Here, ``$q$ is free for $p$ in $\varphi$'' means that no free occurrence of $p$ in $\varphi$ lies within the scope of $\all q$ or $\some q$ (where $\some p$ abbreviates $\neg\all p\neg$).}

\[
\begin{array}{ll}
\text{PC} & \text{all instances of propositional tautologies}\\
\text{$\forall$E} & \all p\varphi\to\varphi(q/p),\text{ where }q\text{ is a propositional variable free for }p\text{ in }\varphi\\
\text{MP} & \dfrac{\varphi\quad \varphi\to\psi}{\psi}\\
\text{$\forall$I} & \dfrac{\varphi\to\psi}{\varphi\to\all p\psi},\ \text{provided }p\text{ does not occur free in }\varphi\\
\text{Mono.} & \dfrac{\varphi\to\psi}{\K{a}\varphi\to\K{a}\psi}
\end{array}
\]

The notions of provability and derivations are defined as usual. We use $\vdash\varphi$ to mean that $\varphi$ is provable in this system.

Note that our system is very weak. For instance, the axiom $\forall$E only permits instantiation by propositional variables free for the quantified variable, rather than arbitrary formulas. By contrast, in Fine's system $\text{AxS5}\pi$ \cite{fine1970}, universal instantiation permits any formula free for the quantified variable. The reason for this restriction is to block the validity of $\all p\varphi\to\varphi(p\land q/p)$ and $\all p\varphi\to\varphi(p\lor q/p)$, which would close the domain of propositions under conjunction and disjunction. This restriction is not ad hoc: it reflects a principled conception of propositional quantification, on which the domain need not be closed under Boolean operations. We provide the philosophical motivation in Section~\ref{sec.conclusion}.

Let $\BP$ be a nonempty set of propositional variables, and $\Ag$ be a nonempty set of agents. The language involved, denoted $\mathcal{L}$, is defined recursively by the following grammar:
\[
\varphi::=p\mid\neg\varphi\mid(\varphi\land\varphi)\mid\K{a}\varphi\mid\all p\varphi,
\]
where $p\in\BP$, $a\in\Ag$.

The above language is interpreted over structures. Intuitively, since a proposition is identified with the set of states at which it is true, propositional quantifiers range over a designated collection of subsets of $S$. Formally, a \emph{structure} $\M$ is a tuple $\langle S,P,R,V\rangle$, where $S$ is a nonempty set of states; $P$ is a nonempty set of subsets of $S$, called the \emph{domain of propositions}; for each $a\in\Ag$, $R_a$ is a binary relation on $S$, called the \emph{accessibility relation}; and $V:\BP\to P$ is a valuation. Truth is defined recursively as follows. Given a structure $\M=\langle S,P,R,V\rangle$ and a state $s\in S$,
\[
\begin{array}{lll}
\M,s\forces p & \text{iff} & s\in V(p)\\
\M,s\forces\neg\varphi & \text{iff} & \M,s\nforces\varphi\\
\M,s\forces\varphi\land\psi & \text{iff} & \M,s\forces\varphi\text{ and }\M,s\forces\psi\\
\M,s\forces\K{a}\varphi & \text{iff} & \text{for all }t,\text{ if }sR_at\text{, then }\M,t\forces\varphi\\
\M,s\forces\all p\varphi & \text{iff} & \text{for all structures }\M'=\langle S,P,R,V'\rangle\\
& & \text{such that }V'(q)=V(q)\text{ for all }q\neq p,\\
& & \M',s\forces\varphi.
\end{array}
\]

Normally, $P$ is subject to various conditions: for instance, that $P$ is Boolean, closed under formulas, or the power set of $S$ \cite[p. 337]{fine1970}. Different conditions may give rise to different logics. For our purposes, we impose no restrictions on $P$. In particular, we require neither closure under intersection nor closure under union. The only requirement is that $V'(p)\in P$ for every valuation $V'$ and $p\in\BP$, which is guaranteed by the definition of structures. This will be crucial in the proofs below. The system $\mathbb{S}$ is sound with respect to this semantics (the verification is routine).

With this machinery in place, we are now in a position to analyze Humberstone's Paradox and its variations in a weak propositionally quantified modal logic, in which fewer theorems are provable since the absence of restrictions on $P$ yields many potential countermodels. We start with Humberstone's Paradox, namely the argument from (1) to (6). Following Humberstone, we reformulate the argument as follows, as is done by Updike~\cite[p. 1185]{updike2024}.

\begin{enumerate}[(i)]
\item $\all p(p\to(\K{a}p\lor\K{b}p))$ \hfill Assumption
\item $(p\land q)\to(\K{a}(p\land q)\lor\K{b}(p\land q))$ \hfill (i), $\forall$E
\item $\K{a}(p\land q)\to\K{a}p$ \hfill by Mono. applied to $(p\land q)\to p$
\item $\K{b}(p\land q)\to\K{b}q$ \hfill by Mono. applied to $(p\land q)\to q$
\item $(p\land q)\to(\K{a}p\lor\K{b}q)$ \hfill (ii),(iii),(iv), PC
\item $(p\to\K{a}p)\lor(q\to\K{b}q)$ \hfill (v), PC
\item $\all q((p\to\K{a}p)\lor(q\to\K{b}q))$ \hfill (vi), $\forall$I
\item $\all p\all q((p\to\K{a}p)\lor(q\to\K{b}q))$ \hfill (vii), $\forall$I
\item $\all p(p\to\K{a}p)\lor\all q(q\to\K{b}q)$ \hfill (viii), Prenex Normal Form Theorem
\end{enumerate}

Note that the above argument is valid if $P$ is closed under intersection.\footnote{The step from (i) to (ii) would be semantically justified if $P$ were closed under intersection, since then the truth-set of $p\land q$ would be in $P$, making $p\land q$ an admissible substitution for $\forall$E.\label{fn.closed-under-intersection}} Without this closure condition, we cannot obtain (ii) from (i) directly. Rather, using $\forall$E with propositional variables only, (i) yields $p\to(\K{a}p\lor\K{b}p)$ and $q\to(\K{a}q\lor\K{b}q)$. We then get $(p\land q)\to(\K{a}p\lor\K{b}p)\land(\K{a}q\lor\K{b}q)$. Deriving (ii) would require $(\K{a}p\lor\K{b}p)\land(\K{a}q\lor\K{b}q)\to\K{a}(p\land q)\lor\K{b}(p\land q)$. But this implication is not valid, as shown by the pointed model below: $\K{a}p$ and $\K{b}q$ are true while $\K{b}p$ and $\K{a}q$ are false, making the antecedent true and the succedent false. Therefore the step from (i) to (ii) is blocked.

Consider the model $\M=\langle S,P,R,V\rangle$, where $S=\{s,t,u\}$, $R_a=\{(s,t)\}$, $R_b=\{(s,u)\}$, $P=\{\{s,t\},\{s,u\}\}$, $V(p)=\{s,t\}$, and $V(q)=\{s,u\}$.  One may check that in $(\M,s)$, (i) is true but (ii) is false. To see that (i) holds, note that for any $V'$ with $V'(p)\in P$, $V'(p)$ is either $\{s,t\}$ or $\{s,u\}$. In either case, $\M',s\forces p\to(\K{a}p\lor\K{b}p)$ (note that since $R_a(s)=\{t\}$ and $R_b(s)=\{u\}$, if $V'(p)=\{s,t\}$ then $t\in V'(p)$, so $\M',s\forces\K{a}p$; if $V'(p)=\{s,u\}$ then $u\in V'(p)$, so $\M',s\forces\K{b}p$). To see that (ii) fails, note that in $(\M,s)$, $p\land q$ is true but $\K{a}(p\land q)$ is false (since $q$ is false at $t$) and $\K{b}(p\land q)$ is false (since $p$ is false at $u$). This shows that (i) does not imply (ii). Furthermore, (ix) is false in $(\M,s)$. If (ix) were true, then either $\M,s\forces\all p(p\to\K{a}p)$ or $\M,s\forces\all q(q\to\K{b}q)$. If the former holds, then $\langle S,P,R,V_1\rangle,s\forces p\to\K{a}p$ for $V_1(p)=\{s,u\}$; since $s\in V_1(p)$, it follows that $\langle S,P,R,V_1\rangle,t\forces p$, hence $t\in V_1(p)$, a contradiction. If the latter holds, then $\langle S,P,R,V_2\rangle,s\forces q\to\K{b}q$ for $V_2(q)=\{s,t\}$; since $s\in V_2(q)$, it follows that $\langle S,P,R,V_2\rangle,u\forces q$, hence $u\in V_2(q)$, again a contradiction.

We now turn to the second argument, from (1$'$) to (6$'$). The proof system is the same as above, with $\B{}$ in place of $\K{}$. As before, we reformulate it with explicit propositional quantifiers:

\begin{enumerate}[(i$'$)]
\item $\all p((\B{a}p\land\B{b}p)\to p)$ \hfill Assumption
\item $(\B{a}(p\lor q)\land\B{b}(p\lor q))\to(p\lor q)$ \hfill (i$'$), $\forall$E
\item $\B{a}p\to\B{a}(p\lor q)$ \hfill by Mono. applied to $p\to(p\lor q)$
\item $\B{b}q\to\B{b}(p\lor q)$ \hfill by Mono. applied to $q\to(p\lor q)$
\item $(\B{a}p\land\B{b}q)\to(p\lor q)$ \hfill (ii$'$), (iii$'$),(iv$'$), PC
\item $(\B{a}p\to p)\lor(\B{b}q\to q)$ \hfill (v$'$), PC
\item $\all q((\B{a}p\to p)\lor(\B{b}q\to q))$ \hfill (vi$'$), $\forall$I
\item $\all p\all q((\B{a}p\to p)\lor(\B{b}q\to q))$ \hfill $(vii')$, $\forall$I
\item $\all p(\B{a}p\to p)\lor\all q(\B{b}q\to q)$ \hfill (viii$'$), Prenex Normal Form Theorem
\end{enumerate}

The above argument is valid if $P$ is closed under union.\footnote{The reason is similar to that in note~\ref{fn.closed-under-intersection}. In detail, the step from (i$'$) to (ii$'$) would be semantically justified if $P$ were closed under union, since then the truth-set of $p\lor q$ would be in $P$, making $p\lor q$ an admissible substitution for $\forall$E.} Without this closure condition, we cannot obtain (ii$'$) from (i$'$) directly. Rather, using $\forall$E with propositional variables only, (i$'$) yields $(\B{a}p\land\B{b}p)\to p$ and $(\B{a}q\land\B{b}q)\to q$. We then get $((\B{a}p\land\B{b}p)\lor(\B{a}q\land\B{b}q))\to(p\lor q)$. Deriving (ii$'$) would require $(\B{a}(p\lor q)\land\B{b}(p\lor q))\to(\B{a}p\land\B{b}p)\lor(\B{a}q\land\B{b}q)$. But this implication is not valid, as shown by the pointed model below: $\B{a}p$ and $\B{b}q$ are true while $\B{b}p$ and $\B{a}q$ are false, making the antecedent true and the succedent false. Therefore the step from (i$'$) to (ii$'$) is blocked.

Consider the model $\M=\langle S,P,R,V\rangle$, where $S=\{s,t,u\}$, $R_a=\{(s,t)\}$, $R_b=\{(s,u)\}$, $P=\{\{t\},\{u\}\}$, $V(p)=\{t\}$, and $V(q)=\{u\}$. One may check that in $(\M,s)$, (i$'$) is true but (ii$'$) is false. To see that (i$'$) holds, note that for any $V'$ with $V'(p)\in P$, $V'(p)$ is either $\{t\}$ or $\{u\}$. In either case, $\M',s\forces(\B{a}p\land\B{b}p)\to p$ (note that if $V'(p)=\{t\}$, then $\B{b}p$ is false at $s$ since $u\notin V'(p)$, so $\B{a}p\land\B{b}p$ is false; if $V'(p)=\{u\}$, then $\B{a}p$ is false at $s$ since $t\notin V'(p)$, so $\B{a}p\land\B{b}p$ is false). To see that (ii$'$) fails, note that in $(\M,s)$, $\B{a}(p\lor q)$ is true (since $p$ is true at $t$) and $\B{b}(p\lor q)$ is true (since $q$ is true at $u$) but $p$ and $q$ are false at $s$. This establishes that (i$'$) does not imply (ii$'$). Moreover, (ix$'$) is false in $(\M,s)$. If (ix$'$) were true, then either $\M,s\forces\all p(\B{a}p\to p)$ or $\M,s\forces\all q(\B{b}q\to q)$. If the former holds, then $\langle S,P,R,V_1\rangle,s\forces\B{a}p\to p$ for $V_1(p)=\{t\}$; since $R_a(s)=\{t\}$ and $t\in V_1(p)$, we have $\langle S,P,R,V_1\rangle,s\forces\B{a}p$, and hence $\langle S,P,R,V_1\rangle,s\forces p$, so $s\in V_1(p)$, a contradiction. If the latter holds, then $\langle S,P,R,V_2\rangle,s\forces\B{b}q\to q$ for $V_2(q)=\{u\}$; since $R_b(s)=\{u\}$ and $u\in V_2(q)$, we have $\langle S,P,R,V_2\rangle,s\forces\B{b}q$, and hence $\langle S,P,R,V_2\rangle,s\forces q$, so $s\in V_2(q)$, again a contradiction.

We have therefore demonstrated that, in the weak propositionally quantified modal logic $\mathbb{S}$ (which is consistent), Humberstone's Paradox and its variations are blocked: for each paradox we have constructed a (pointed) model in which the premise holds but the conclusion fails. This block extends also to the distribution-free variant of Humberstone's Paradox developed by Updike \cite{updike2024}, since both key lemmas of his reconstruction require compound instantiation of universally quantified premises (collective omniscience in Lemma 1, the knowability principle in Lemma 2), and such instantiation is respectively unavailable under our restricted \(\forall\)E and \(\forall\)E\('\) axioms.

\medskip

Now we move to the analysis of Fitch's Paradox. The situation now becomes complicated, since here the required system is stronger and the semantics more complex. In detail, let $\mathbb{S}'$ consist of the following axioms and inference rules, where the class of \emph{special forms} is defined recursively: every propositional variable is a special form, and if $\varphi$ is a special form, then so are $\neg\varphi$ and $\K{}\varphi$. Thus special forms are exactly the formulas built from propositional variables using only $\neg$ and $\K{}$; they contain no $\land$, $\lor$, $\Boxm$, or $\Dia$. This class is semantically motivated: since $P$ is closed under relative complement and $m_K$ but not under $m_\Box$, intersection, or union, the truth-set of every special form is guaranteed to lie in $P$, whereas the truth-set of a formula containing $\land$, $\lor$, or $\Boxm$ need not.

\[
\begin{array}{ll}
\text{PC} & \text{all instances of propositional tautologies}\\
\forall\text{E}' & \all p\varphi\to\varphi(\psi/p),\text{ where }\psi\text{ is a special form free for }p\text{ in }\varphi\\
\text{T} & \K{}\varphi\to\varphi\\
\text{MP} & \dfrac{\varphi\quad \varphi\to\psi}{\psi}\\
\text{$\forall$I} & \dfrac{\varphi\to\psi}{\varphi\to\all p\psi},\ \text{provided }p\text{ does not occur free in }\varphi\\
\text{RN} & \dfrac{\varphi}{\Boxm\varphi}\\
\text{RE} & \dfrac{\varphi\leftrightarrow\psi}{\Boxm\varphi\leftrightarrow\Boxm\psi}\\
\text{Mono.} & \dfrac{\varphi\to\psi}{\K{}\varphi\to\K{}\psi}
\end{array}
\]

The language involved, denoted $\mathcal{L}'$, is an extension of $\mathcal{L}$ with the additional operator $\Boxm$. In detail, $\mathcal{L}'$ is defined recursively by the following grammar:
\[
\varphi::=p\mid\neg\varphi\mid(\varphi\land\varphi)\mid\K{}\varphi\mid\all p\varphi\mid\Boxm\varphi,
\]
where $p\in\BP$. Other connectives are defined as usual; in particular, $\Dia\varphi$ abbreviates $\neg\Boxm\neg\varphi$.

Semantically, $\mathcal{L}'$ is interpreted over structures. A structure $\M$ is a tuple $\langle S,R_\Boxm,R_K,P,V\rangle$, where $S,P,V$ are defined as before, $R_\Boxm$ and $R_K$ are binary relations on $S$, and $R_K$ is reflexive.\footnote{Recall that in Fitch's original argument, knowledge is factive, which corresponds to the reflexivity of $R_K$.} Given a structure $\M=\langle S,R_\Boxm,R_K,P,V\rangle$ and a state $s\in S$, the semantics of $\mathcal{L}'$ extends that of $\mathcal{L}$ with
\[
\begin{array}{lll}
\M,s\forces\Boxm\varphi & \Longleftrightarrow & \text{for all }t,\text{ if }sR_\Boxm t\text{ then }\M,t\forces\varphi\\
\M,s\forces\K{}\varphi & \Longleftrightarrow & \text{for all }t,\text{ if }sR_K t\text{ then }\M,t\forces\varphi
\end{array}
\]

The semantics of $\Dia\varphi$ follows immediately:
\[
\begin{array}{lll}
\M,s\forces\Dia\varphi & \Longleftrightarrow & \text{there exists }t\text{ such that }sR_\Boxm t\text{ and }\M,t\forces\varphi
\end{array}
\]

For our purposes of resolving Fitch's Paradox, we require that $P$ satisfy the following conditions:
\begin{itemize}
\item closed under relative complement: if $X\in P$, then $S\setminus X\in P$.
\item closed under modal operation $m_K$: if $X\in P$, then $m_K(X)\in P$, where $m_K(X)=\{s\in S\mid\text{for all }t,\text{ if }sR_K t\text{ then }t\in X\}$.
\end{itemize}
Analogously, define $m_\Box(X)=\{s\in S\mid\text{for all }t,\text{ if }sR_\Boxm t\text{ then }t\in X\}$. We do not require $P$ to be closed under $m_\Box$.

It is then straightforward to verify that $V(\K{}\varphi)=m_K(V(\varphi))$. We now reconstruct Fitch's Paradox as follows:
\begin{enumerate}[(i$''$)]
\item $\all p(p\to\Dia\K{}p)$ \hfill Assumption
\item $(p\land\neg\K{}p)\to\Dia\K{}(p\land\neg\K{}p)$ \hfill (i$''$), $\forall$E$'$
\item $\K{}(p\land\neg\K{}p)\to(\K{}p\land\K{}\neg\K{}p)$ \hfill Mono., PC
\item $\K{}\neg\K{}p\to\neg\K{}p$ \hfill Ax.~T
\item $\K{}(p\land\neg\K{}p)\to(\K{}p\land\neg\K{}p)$ \hfill (iii$''$), (iv$''$), PC
\item $\neg(\K{}p\land\neg\K{}p)$ \hfill PC
\item $\neg\K{}(p\land\neg\K{}p)$ \hfill (v$''$), (vi$''$), PC
\item $\Boxm\neg\K{}(p\land\neg\K{}p)$ \hfill (vii$''$), RN
\item $\neg\Dia\K{}(p\land\neg\K{}p)$ \hfill (viii$''$), Def.~$\Dia$, RE
\item $\neg(p\land\neg\K{}p)$ \hfill (ix$''$), (ii$''$), PC
\item $p\to\K{}p$ \hfill (x$''$), PC
\item $\all p(p\to\K{}p)$ \hfill (xi$''$), $\forall$I
\end{enumerate}

The above argument is valid if $P$ is closed under intersection.\footnote{The reason is analogous to the Humberstone case. In detail, the step from (i$''$) to (ii$''$) would be semantically justified if $P$ were closed under intersection, since then the truth-set of $p\land\neg\K{}p$ would be in $P$, making $p\land\neg\K{}p$ an admissible substitution for $\forall$E$'$.} Without this closure condition, we cannot obtain (ii$''$) from (i$''$) directly. Rather, using $\forall$E$'$ with special forms only, $(i'')$ yields $p\to\Dia\K{}p$ and $\neg\K{}p\to\Dia\K{}\neg\K{}p$. We then get $(p\land\neg\K{}p)\to(\Dia\K{}p\land\Dia\K{}\neg\K{}p)$. Deriving $(ii'')$ would require $\Dia\K{}p\land\Dia\K{}\neg\K{}p\to\Dia\K{}(p\land\neg\K{}p)$. But this implication is not valid, as shown by the pointed model below, where the antecedent is true but the succedent is false. Therefore the step from (i$''$) to (ii$''$) is blocked.

Consider the model $\M=\langle S,R_\Boxm,R_K,P,V\rangle$ where $S=\{s,t,u\}$, $R_\Boxm=\{(s,t),(s,u)\}$, $R_K=\{(s,s),(t,t),(u,u),(s,t),(s,u)\}$, $P=\{\{s,t\},\{u\},\{t\},\{s,u\}\}$, $V(p)=\{s,t\}$. It is straightforward to verify that $P$ is closed under relative complement and modal operation $m_K$. Moreover, in $(\M,s)$, (i$''$) is true but (ii$''$) is false. To see that (i$''$) holds, note that for any $V'$ with $V'(p)\in P$, $V'(p)$ is one of $\{s,t\},\{u\},\{t\},\{s,u\}$. In each case, $\M',s\forces p\to\Dia\K{}p$ (if $V'(p)=\{u\}$ or $V'(p)=\{t\}$, then $s\notin V'(p)$, so $\M',s\forces p\to\Dia\K{}p$; if $V'(p)=\{s,t\}$, then $sR_\Boxm t$, $R_K(t)=\{t\}$, and $t\in V'(p)$, so $\M',s\forces\Dia\K{}p$, hence $\M',s\forces p\to\Dia\K{}p$; if $V'(p)=\{s,u\}$, then $sR_\Boxm u$, $R_K(u)=\{u\}$, and $u\in V'(p)$, so $\M',s\forces\Dia\K{}p$, hence $\M',s\forces p\to\Dia\K{}p$). To see that (ii$''$) fails, note that $\M,s\forces p\land\neg\K{}p$ and $\M,s\nforces\Dia\K{}(p\land\neg\K{}p)$ (otherwise, $t\forces\K{}(p\land\neg\K{}p)$ or $u\forces\K{}(p\land\neg\K{}p)$, either of which yields a contradiction). This establishes that (i$''$) does not imply (ii$''$). Moreover, from $\M,s\forces p\land\neg\K{}p$, it follows that $\M,s\nforces p\to\K{}p$, and hence $\M,s\nforces\all p(p\to\K{}p)$, so (xii$''$) is also false in $(\M,s)$.

Our approach to resolving Fitch's Paradox is similar to that of Restall \cite{restall2009}: even if both $p$ and $\neg\K{}p$ are knowable, $p\land\neg\K{}p$ is unknowable. There are, however, important differences. Restall \cite{restall2009} does not endorse the knowability thesis, that is, (i$''$). Rather, he endorses what he calls the ``conjunctive knowability thesis'', according to which every truth is conjunctively knowable; that is, it is equivalent to a conjunction each of whose conjuncts is knowable.  This is weaker than the knowability thesis. In contrast, we endorse the knowability thesis. Moreover, Restall accepts the conclusion of Fitch's Paradox, whereas we reject it. We do not merely block the argument: we have constructed a countermodel showing that the knowability thesis holds while omniscience fails.

Another related approach is Dummett's \cite{dummett2001} (see also \cite[p. 398]{humberstone2015book}), where he proposes restricting the knowability principle to basic (or atomic) statements. He gives an inductive characterization of truth, on which the knowability principle applies only to basic statements, not, for example, to conjunctions thereof. Thus he endorses, for instance, $p\to\Dia\K{}p$ and $q\to\Dia\K{}q$, but not $p\land q\to\Dia\K{}(p\land q)$. This syntactic restriction allows Dummett to avoid Fitch's Paradox.

Our approach differs from Dummett's in a crucial respect. Dummett restricts the knowability principle to basic statements, whereas we retain the principle in its full generality as $\all p(p\to\Dia\K{}p)$ and instead weaken the propositional quantifier: its instantiation axiom $\forall$E$'$ admits only special forms as substituends. Moreover, Dummett's restriction is syntactic, grounded in an inductive characterization of truth, whereas ours is semantic, grounded in the closure conditions of the domain of propositions. Our framework also provides a uniform resolution of both Fitch's and Humberstone's Paradoxes, whereas Dummett's proposal targets Fitch's Paradox alone.

\section{Conclusion and Responses to Objections}\label{sec.conclusion}
\subsection{Addressing the Motivation and Nature of the Solution}
Before concluding, we address two central concerns regarding the philosophical motivation of our proposal and the nature of our ``resolution'' of the paradoxes.

First, the worry that our framework is \emph{ad hoc} rests on the unstated assumption that the only legitimate semantics for propositional quantifiers requires the domain $P$ to be closed under all Boolean operations and modal constructors. This assumption is not forced by the standard account of propositional quantification: as Fine \cite{fine1970} notes, conditions on $P$ are a matter of substantive philosophical choice, not a purely logical default. Our restriction on universal instantiation ($\forall$E) is not introduced merely to block the paradoxical inferences, but derives from an independently motivated conception of propositions as graspable objects of epistemic attitudes.

To clarify this proposal, we distinguish graspability from knowability. For a single agent with finite cognitive capacities, a proposition is \emph{graspable} if it can serve as the unitary object of a cognitive attitude apt for knowledge for that agent. Graspability requires not merely intelligibility, but cognitive stability: a graspable proposition must be capable of being held as a coherent candidate for knowledge without self-defeat. Thus Moorean propositions such as $p\land\neg\K{}p$ are intelligible but not graspable: although one can understand what they say, they cannot be coherently held as candidates for knowledge, and hence are not knowable, since knowing them would be self-defeating. Graspability is conceptually prior to knowability: a proposition must be graspable to be knowable, but graspability does not guarantee knowability. For example, a proposition about a particular past event that left no recoverable trace is intelligible and cognitively stable, hence graspable, yet not knowable, since no finite agent can in principle gather sufficient evidence to ascertain its truth. We propose that the domain of propositions should contain exactly the graspable propositions, i.e., those that are cognitively stable and eligible as candidates for knowledge. Note that this claim concerns the domain over which epistemic principles such as the knowability thesis quantify, not the broader question of whether every intelligible declarative sentence expresses a proposition; we remain neutral on the latter question. Our point is that a content that cannot in principle be grasped as a candidate for knowledge falls outside the scope of the anti-realist's ``every truth'' quantifier, regardless of whether it qualifies as a proposition in some broader sense. The rationale is that epistemic principles such as the knowability thesis concern what finite agents can come to know, and so should quantify only over propositions that finite agents can grasp as objects of cognitive attitudes; this dovetails with the anti-realist conception of propositions as possible objects of verification or knowledge, on which a proposition that cannot in principle be grasped as a candidate for knowledge is not the sort of thing over which the anti-realist's ``every truth'' quantifies. Because graspability is constrained by the cognitive limits of finite agents rather than by purely logical structure, the set of graspable propositions need not exhibit the Boolean closure properties standardly assumed for propositional domains.

These cognitive constraints also determine which closure properties the domain of graspable propositions does and does not possess. On the one hand, the domain is plausibly closed under negation: if a proposition is graspable, its negation is graspable as well, since negation is a truth-functional operation that requires no additional cognitive resources beyond understanding the original proposition. It is also plausibly closed under the knowledge operation: if a proposition is graspable, then the proposition that it is known is also graspable, since a cognitively competent agent can in principle reflect on whether it knows a given proposition. This closure under $m_K$ should not be confused with endorsement of the KK principle ($\K{}p\to\K{}\K{}p$): it concerns only the graspability of knowledge ascriptions, not their truth. On the other hand, the domain need not be closed under the necessity operation $m_\Box$, since grasping that a proposition is metaphysically necessary requires ontological reasoning that may exceed ordinary cognitive capacities. Nor need it be closed under conjunction or disjunction: unlike negation, which is a unary operation requiring only that one take a proposition already grasped and reverse its truth-value, conjunction and disjunction are binary operations that require simultaneously entertaining two propositions and grasping their relation, which may either exceed the capacity of a finite agent (as in the Humberstone case) or yield a self-defeating content (as in the Fitch case). The asymmetry thus reflects the cognitive difference between unary and binary operations, not an arbitrary stipulation, for the reasons developed below.

The failure of closure arises from two distinct sources, corresponding to the two paradoxes.

\paragraph{The Fitch case: cognitive self-defeat.} In the Fitch case, $p\land\neg\K{}p$ is a cognitive analogue of a Moore sentence: it may be true, but no minimally reflective agent can entertain it as a unitary object of knowledge. Here ``graspable'' is understood in a specific sense: it means eligible as a candidate for knowledge, not merely entertainable or believable, and not synonymous with knowable. As noted above, graspability is conceptually prior to knowability: a proposition must be graspable to be knowable, but graspability does not guarantee knowability (witness propositions about past events that left no recoverable trace). One can certainly entertain, suppose, or even believe $p\land\neg\K{}p$, but one cannot coherently hold it as a candidate for knowledge. The point is not merely that $p\land\neg\K{}p$ is unknowable (a claim that would risk circularity, since unknowability is what Fitch's argument establishes). Rather, the point is that $p\land\neg\K{}p$ cannot be stably grasped as a candidate for knowledge at all. To see this, suppose for \emph{reductio} that $p\land\neg\K{}p$ is known:

\begin{enumerate}[(1)]
\item $\K{}(p\land\neg\K{}p)$ \hfill Assumption
\item $\K{}p\land\K{}\neg\K{}p$ \hfill Knowledge distributes over conjunction
\item $\K{}\neg\K{}p\to\neg\K{}p$ \hfill Factivity of knowledge
\item $\K{}p\land\neg\K{}p$ \hfill A truth-functional consequence of (2) and (3)
\end{enumerate}
Thus the assumption that $p\land\neg\K{}p$ is known leads to a logical contradiction. This is not a merely pragmatic oddity; it is a structural feature of the proposition: its truth condition itself undermines the possibility of its being known. For an anti-realist who ties the notion of truth to the possibility of verification or knowledge, a proposition whose truth would render it unknowable is not a legitimate \emph{truth-candidate}; it does not fall within the range of propositions over which the anti-realist's ``every truth'' quantifies. This is not a stipulation introduced to avoid the paradox, but a consequence of the anti-realist conception of what propositions are: propositions are possible objects of verification, not arbitrary truth-functional constructions.

This approach is in the spirit of Tennant's \emph{Cartesian restriction} on knowability \cite{tennant1997}, according to which the knowability principle applies only to propositions whose being known is consistent (``Cartesian'' propositions). $p\land\neg\K{}p$ is non-Cartesian, since assuming it known yields a contradiction. Where we go beyond Tennant is that we provide this restriction with a \emph{semantic foundation}: rather than screening propositions one by one for Cartesian consistency via a syntactic test, our restricted-domain semantics excludes self-defeating conjunctions globally through the structure of the propositional domain. The restricted $\forall$E$'$ axiom, which permits substitution only of special forms (formulas built from propositional variables using $\neg$ and $\K{}$), directly tracks which formulas correspond to graspable propositions, providing a precise syntactic correlate of our semantic constraint on $P$.

A further source of the \emph{ad hoc} worry concerns the specific scope of this instantiation axiom. One might ask why the special forms are exactly the $\neg,\K{}$-combinations, rather than a broader class that also includes $\Boxm\varphi$, or a narrower class containing only propositional variables. The answer is that the class is determined by the closure conditions on $P$: since $P$ is closed under relative complement and $m_K$, the truth-set of every $\neg,\K{}$-combination is guaranteed to lie in $P$; since $P$ is not closed under $m_\Box$, intersection, or union, the truth-set of a formula containing $\Boxm$, $\land$, or $\lor$ need not lie in $P$. The special forms are thus precisely the formulas whose admissibility as substituends is semantically guaranteed. This instantiation axiom is not an arbitrary syntactic cutoff; it is calibrated to the closure properties of the propositional domain, in that its scope reflects exactly which formulas are guaranteed to express propositions in $P$.

\paragraph{The Humberstone case: finite cognitive resources.} In the Humberstone case, the failure of closure under conjunction has a different source. Grasping a proposition consumes finite cognitive resources: attention, memory, computational capacity, and even spatiotemporal location, as emphasized by the shared-resource intuition underlying Humberstone's Paradox \cite[pp. 1183--1184]{updike2024}. Two propositions may each be graspable by some agent, while their conjunction exceeds the capacity of any single agent. The domain of graspable propositions thus need not be closed under conjunction, not by stipulation, but as a reflection of the structural limits of finite cognition. The same consideration applies to disjunction: grasping a disjunction may require entertaining both disjuncts and their relation, which can likewise exceed the capacity of any single agent. In both cases, then, the failure of Boolean closure has independent philosophical roots rather than being a stipulated restriction.

Our approach is uniform: it resolves both Fitch's and Humberstone's paradoxes as well as their variations involving belief, permissibility, and desire, via the single constraint that $P$ need not be closed under intersection or union. This uniformity further confirms that our restriction captures a shared structural flaw in the paradoxical arguments, rather than constituting a targeted patch for isolated problems.

Second, the concern that our framework merely suppresses the expression of the paradoxes rather than solving them misunderstands the status of the solution. We do not deny that one can formulate absolute universal quantification over an unrestricted domain of propositions in other logical systems; our point is that such unrestricted quantification misrepresents the content of the knowability principle and the collective omniscience principle. When these principles are stated informally, as in ``every truth is knowable'' and ``every truth is known by someone'', the quantifier ranges over truths that are, by the very lights of these principles, graspable. Our framework captures precisely this intended reading: the premise \(\all p(p\to\Dia\K{}p)\) holds in our models exactly when every graspable true proposition is knowable, which, we argue, captures the intended content of the anti-realist's knowability principle; the collective omniscience principle \(\all p(p\to\K{a}p\lor\K{b}p)\) likewise holds exactly when every graspable true proposition is known by at least one agent. What our framework blocks is not the statement of the premises, but the invalid inference from these reasonable premises to the absurd conclusions, by showing that the step from a quantified premise to the instantiated conjunctive or disjunctive formula relies on an unwarranted expansion of the propositional domain. This is not a linguistic evasion, but a diagnostic clarification: the paradoxes arise not from the core intuitions behind the knowability or collective omniscience principles, but from an overpermissive conception of propositions tacitly assumed in their standard quantified reconstructions.

We emphasize that we do not take Fitch's or Humberstone's informal arguments to be fallacious in all contexts. They are perfectly coherent as expressions of the relevant epistemic commitments. The problem arises only when these commitments are formalized using propositional quantifiers with an overpermissive domain closure assumption. Our proposal thus offers not a rejection of these commitments, but a more accurate formalization of them, one that preserves their intuitive force while avoiding the absurd consequences. This qualifies as a genuine resolution of the paradoxes: we identify the hidden presupposition that generates the contradiction, show that this presupposition is philosophically dispensable, and provide a unified, semantically grounded alternative framework in which the original motivating principles remain tenable.

In sum, our approach advances the literature in two key respects. First, it provides a unified diagnosis of two long-standing epistemic paradoxes, showing they share a common root in unexamined domain-closure assumptions. Second, it grounds the previously \emph{ad hoc} syntactic restrictions for Fitch's paradox in a rigorous model theory, addressing the long-standing worry that such restrictions lack independent motivation. We hope this clarifies that our framework is not a piecemeal logical trick, but a principled rethinking of how propositional quantification ought to be deployed in epistemic contexts.

\subsection{Summary of the Technical Contribution}

To summarize the technical core of our proposal, we have demonstrated that the paradoxical inferences in both Fitch's and Humberstone's arguments hinge crucially on a tacit assumption regarding the domain of propositions, namely that this domain is closed under intersection and union. By constructing deliberately weak propositionally quantified modal logics, we have shown that relinquishing this assumption provides a uniform and effective blockage of the paradoxes.

Syntactically, the crux of our system lies in the restriction on universal instantiation ($\forall$E in the Humberstone case, $\forall$E$'$ in the Fitch case). Unlike standard systems where $\all p\varphi$ entails $\varphi(\psi/p)$ for any well-formed $\psi$, our framework limits substitution to forms that correspond to graspable propositions: propositional variables in the multi-agent Humberstone case, and special forms (formulas built from propositional variables using $\neg$ and $\K{}$) in the Fitch case. This restriction is not merely a technical expedient; it is the precise syntactic correlate of our semantic decision not to require closure of the propositional domain $P$ under intersection or union. Consequently, the step from the general premise (e.g., $\all p(p\to\K{a}p\lor\K{b}p)$) to the specific instantiation involving a conjunction or disjunction (e.g., $p\land q\to\K{a}(p\land q)\lor\K{b}(p\land q)$) is blocked, as the truth-set of such a compound formula is not guaranteed to reside within $P$.

This approach offers two distinct advantages over extant solutions. First, it achieves unification: the same mechanism, denying closure under intersection and union, resolves Humberstone's Paradox, its variations concerning belief, desire, and permissibility, and Fitch's Paradox simultaneously. This suggests that these puzzles are not isolated anomalies but symptoms of a shared over-generation of propositions in standard reconstructions. Second, it provides a semantic foundation for syntactic restrictions: where previous proposals (such as Dummett's) faced charges of being \emph{ad hoc} for limiting the knowability principle to atomic statements, our restricted-domain semantics supplies an independent model-theoretic motivation. The restricted instantiation axiom in our system, while non-standard, is not without precedent: Fine \cite{fine1970} already considers propositionally quantified modal logics with varying strengths of the universal instantiation axiom, the broader literature on propositional quantification (e.g. \cite{fritz2024}) explores a range of domain conditions, and failures of uniform substitution also arise in dynamic epistemic logic, illustrating that non-standard substitution behavior is not unique to propositionally quantified systems \cite{vanditmarsch2007}. This restriction serves to explain why the paradoxes arise in the first place.

Crucially, our analysis does not impugn the informal validity of Fitch's or Humberstone's original, unquantified reasoning \emph{per se}. Rather, we diagnose a mismatch between their intuitive epistemic commitments and the standard formal tools used to regiment them. The implicit quantification over an unrestricted domain of propositions imports ontological commitments, specifically closure under intersection and union, that were absent from the original informal principles. By refining the underlying logic to respect the constraints on cognitive grasp, we preserve the spirit of anti-realism and collective epistemology without succumbing to omniscience. In this sense, our work offers not merely a technical fix, but a clarification of the logical architecture required for a coherent treatment of propositional quantification in epistemic contexts.

\section*{Acknowledgements}

This article is supported by the Key Project of National Social Science Fund of China (No. 24AZX019). The author thanks Lloyd Humberstone for his encouragement and for many helpful comments on earlier drafts of this paper, and the anonymous referee of this journal for insightful and constructive comments on the submitted version.


\begin{thebibliography}{99}

\bibitem{bigelow2005} J. Bigelow. Omnificence. \textit{Analysis}, 65(3):187--196, 2005.

\bibitem{brogaard2019} B. Brogaard and J. Salerno. Fitch's paradox of knowability. In Edward N. Zalta, editor, \textit{The Stanford Encyclopedia of Philosophy}. Metaphysics Research Lab, Stanford University, fall 2019 edition, 2019.

\bibitem{dummett2001} M. Dummett. Victor's error. \textit{Analysis}, 61:1--2, 2001.

\bibitem{dummett2009} M. Dummett. Fitch's paradox of knowability. In J. Salerno, editor, \textit{New essays on the knowability paradox}, pages 51--52. Oxford University Press, 2009.

\bibitem{fine1970} K. Fine. Propositional quantifiers in modal logic. \textit{Theoria}, 36(3):336--346, 1970.

\bibitem{fitch1963} F. B. Fitch. A logical analysis of some value concepts. \textit{The Journal of Symbolic Logic}, 28(2):135--142, 1963.

\bibitem{fritz2024} P. Fritz. \textit{Propositional Quantifiers}. Cambridge University Press, 2024.


\bibitem{guigon2009} G. Guigon. Bringing about and conjunction: A reply to Bigelow on omnificence. \textit{Analysis}, 69(3):452--458, 2009.

\bibitem{humberstone1985} L. Humberstone. The formalities of collective omniscience. \textit{Philosophical Studies}, 48(3):401--423, 1985.

\bibitem{humberstone2015book} L. Humberstone. \textit{Philosophical Applications of Modal Logic}. London: College Publications, 2015.

\bibitem{maglo2010} M. Jago. Closure on knowability. \textit{Analysis}, 70(4):648--659, 2010.

\bibitem{linsky1986} B. Linsky. Factives, blindspots and some paradoxes. \textit{Analysis}, 64:10--15, 1986.

\bibitem{restall2009} G. Restall. Not every truth can be known (at least, not all at once). In J. Salerno, editor, \textit{New essays on the knowability paradox}, pages 339--354. Oxford University Press, 2009.

\bibitem{salerno2009} J. Salerno. Knowability noir: 1945--1963. In J. Salerno, editor, \textit{New Essays on the Knowability Paradox}, pages 29--48. Oxford University Press, 2009.

\bibitem{salerno2009book} J. Salerno, editor. \textit{New Essays on the Knowability Paradox}. Oxford University Press, Oxford, UK, 2009.

\bibitem{tennant1997} N. Tennant. \textit{The Taming of the True}. Oxford University Press, Oxford, 1997.

\bibitem{updike2022} E. Updike. Collecting truths: A paradox in two guises. \textit{Analytic Philosophy}, 63(3):156--173, 2022.

\bibitem{updike2024} E. Updike. Humberstone's paradox and conjunction. \textit{Erkenntnis}, 89:1183--1195, 2024.

\bibitem{vanbenthem2009} J. van Benthem. Actions that make us know. In J. Salerno, editor, \textit{New essays on the knowability paradox}, pages 129--146. Oxford University Press, 2009.

\bibitem{vanditmarsch2007} H. van Ditmarsch, W. van der Hoek, and B. Kooi. 
\textit{Dynamic Epistemic Logic}, volume 337 of Synthese Library. Springer, 2007.


\bibitem{williamson1993} T. Williamson. Verificationism and non-distributive knowledge. \textit{Australasian Journal of Philosophy}, 71:78--86, 1993.

\end{thebibliography}
\end{document}